%% file: paper.tex
\documentclass[final]{siamart220329}

\usepackage{lipsum}
\usepackage{amsfonts}
\usepackage{graphicx}
\usepackage{epstopdf}
\usepackage{algorithmic}
\ifpdf%
  \DeclareGraphicsExtensions{.eps,.pdf,.png,.jpg}
\else
  \DeclareGraphicsExtensions{.eps}
\fi
\usepackage{amsopn}
\DeclareMathOperator{\diag}{diag}
\usepackage{booktabs}
\usepackage[mathscr]{eucal}
\input{tikz_tensor_object.tex}

\newcommand{\T}[1]{\boldsymbol{\mathscr{\MakeUppercase{#1}}}}
\newcommand{\M}[1]{\mathbf{\MakeUppercase{#1}}}
\newcommand{\V}[1]{\mathbf{#1}}
\newcommand{\bmp}{{\mbox{\tt bmp}}}
\newcommand{\tvec}{{\mbox{\tt Tvec}}}
\newcommand{\tfold}{{\mbox{\tt Tfold}}}
\newcommand{\mat}{{\mbox{\tt Mat}}}

\newcommand{\sq}{{\mbox{\tt squeeze}}}
\newcommand{\bdiag}{{\mbox{\tt bdiag}}}
\newcommand{\perm}{{\mbox{\tt permute}}}

\newcommand{\TheTitle}{%
Tensor Completion with BMD Factor Nuclear Norm Minimization
}

\newcommand{\TheShortTitle}{%
  \TheTitle
}

\newcommand{\TheName}{%
  Fan Tian,   Mirjeta Pasha,
  Misha E. Kilmer,
  Eric Miller,
  Abani Patra
}

\newcommand{\TheFunding}{%
This work was partially supported by the National Science
Foundation under NSF HDR grant CCF-1934553. Misha E. Kilmer was also supported in part by NSF DMS-1821148. MP gratefully acknowledges support from the NSF under award No. 2202846.
}

\author{\TheName}
\title{{\TheTitle}\thanks{\TheFunding}}
\headers{\TheShortTitle}{\TheName}
\ifpdf%
\hypersetup{%
  pdftitle={\TheTitle},
  pdfauthor={\TheName}
}
\fi

\begin{document}

\maketitle


\begin{abstract}
This paper is concerned with the problem of recovering third-order tensor data from limited samples. A recently proposed tensor decomposition (BMD) method has been shown to efficiently compress third-order spatiotemporal data. Using the BMD, we formulate a slicewise nuclear norm penalized algorithm to recover a third-order tensor from limited observed samples. We develop an efficient alternating direction method of multipliers (ADMM) scheme to solve the resulting minimization problem. Experimental results on real data show our method to give reconstruction comparable to those of HaLRTC (Liu et al., IEEE Trans Ptrn Anal Mchn Int, 2012), a well-known tensor completion method, in about the same number of iterations. However, our method has the advantage of smaller subproblems and higher parallelizability per iteration.

\end{abstract}

\begin{keywords}
  tensor completion, tensor BM-decomposition, ADMM
\end{keywords}

\section{Introduction}\label{sec:intro}
Missing data can arise in many scenarios in real-world video and image processing applications. For instance, malfunctioning sensors or cameras can lead to missing data in photos and video recordings. Environmental factors like adverse weather conditions, poor lighting, or occlusions can affect the quality of the data and result in missing or distorted frames. Studying higher-order tensor completion problems incorporating the multidimensional structure of the color images, hyperspectral images or grayscale videos has drawn great attention in recent years \cite{liu2012tensor, zhang2014novel}. A survey on the recent developments in methods and algorithms for solving the tensor completion problem is available in \cite{song2019tensor}.

A common assumption for solving the tensor completion problem is that the underlying unknown data tensor can be well approximated by a low rank tensor. The definition of the tensor rank can vary depending on the underlying rank assumption of a given tensor \cite{song2019tensor}. Commonly studied tensor ranks include the tensor CP rank \cite{kolda2009tensor}, the Tucker rank \cite{de2000multilinear}, and the tensor tubal-rank known from the t-SVD factorization \cite{kilmer2011factorization}. 

More recently, a novel tensor Bhattacharya-Mesner decomposition (BMD) has been proposed  \cite{tian2023tensor} for third-order tensors. It is demonstrated that third-order spatiotemporal data is highly compressible. In particular, we observed that the spatiotemporal slices of the surveillance video exhibit low-rankness. Moreover, a generative low BM-rank video model was introduced and has shown that the motions of objects traversing the background scene over time can be captured within a few slices of the BMD factor tensors. We expect that these factor tensor slices are low rank. As for the factor tensor slices corresponding to the spatial dimension, they are less expected to be low rank. Motivated by these observations, we propose a novel BM-factor tensor slicewise nuclear norm minimization (BMNN) method for the tensor completion problem. In this work, we focus on third-order tensors.

This paper is organized as follows. In section \ref{sec:bg}, we present the notations and an overview of the tensor BMD. In section \ref{sec:main}, we formulate the tensor completion problem with tensor BMD and introduce a tractable relaxation of the rank function into our factor tensor slicewise rank minimization model. Then we will provide an efficient alternating direction method of multipliers (ADMM) scheme to solve the proposed problem. In section \ref{sec:numerical}, we will present our experimental findings applied to real-world videos and hyperspectral images. We compare our reconstruction performances with a well-known tensor completion method, namely the HaLRTC algorithm in \cite{liu2012tensor}. Finally, we conclude in section \ref{sec:conclude}.

\section{Background and Notation}\label{sec:bg}
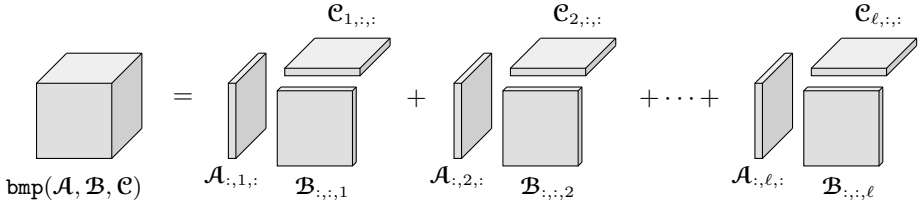
\begin{figure}[H]
\centering
\begin{tikzpicture}
\tensor[dim1=1, dim2=1, dim3=0.8, fill = lightgray] (0.2, 0) {};
\node at (0.7, -0.4) {\small $\bmp(\T{A},\T{B},\T{C})$};
\node at (2.15, 0.8) {$=$};
\tensor[dim1=1, dim2=0.1, dim3=0.8, fill = lightgray] (2.75, 0) {};
\tensor[dim1=1, dim2=1, dim3=0.1, fill = lightgray] (3.4, -0.1) {};
\tensor[dim1=0.1, dim2=1, dim3=0.8, fill = lightgray] (3.5, 1.1) {};
\node at (2.8, -0.25) {\small $\T{\T{A}}_{:,1,:}$};
\node at (4, -0.4) {\small $\T{B}_{:,:,1}$};
\node at (4.4, 1.9) {\small $\T{C}_{1,:,:}$};
\node at (3+2.25, 0.8){$+$};
\tensor[dim1=1, dim2=0.1, dim3=0.8, fill = lightgray] (3+2.75, 0) {};
\tensor[dim1=1, dim2=1, dim3=0.1, fill = lightgray] (3+3.4, -0.1) {};
\tensor[dim1=0.1, dim2=1, dim3=0.8, fill = lightgray] (3+3.5, 1.1) {};
\node at (3+2.8, -0.25) {\small $\T{\T{A}}_{:,2,:}$};
\node at (3+4, -0.4) {\small $\T{B}_{:,:,2}$};
\node at (3+4.4, 1.9) {\small $\T{C}_{2,:,:}$};
\node at (6.6+2.15, 0.8){$+ \cdots +$};
\tensor[dim1=1, dim2=0.1, dim3=0.8, fill = lightgray] (7+2.75, 0) {};
\tensor[dim1=1, dim2=1, dim3=0.1, fill = lightgray] (7+3.4, -0.1) {};
\tensor[dim1=0.1, dim2=1, dim3=0.8, fill = lightgray] (7+3.5, 1.1) {};
\node at (7+2.8, -0.25) {\small $\T{\T{A}}_{:,\ell,:}$};
\node at (7+4, -0.4) {\small $\T{B}_{:,:,\ell}$};
\node at (7+4.4, 1.9) {\small $\T{C}_{\ell,:,:}$};
\end{tikzpicture}
\caption{Illustration of the BM-product of a conformable tensor triplet $\T{A},\T{B}$, and $\T{C}$ as a sum of BM-outer products of matrix slices.}
\label{fig:bmp_outer}
\end{figure}

The tensor BMD is defined based on the BM-product of higher-order tensors, which were first introduced by Bhattacharya and Mesner as a ternary association scheme on triplets \cite{mesner1990association, mesner1994ternary} and later studied by Gnang and Filmus \cite{gnang2017spectra,gnang2020bhattacharya}. For the third-order tensors, the definition of the BM-product is given as follows
\begin{definition}
For a third-order conformable tensor triplet $\T{A} \in \mathbb{R}^{m \times \ell \times p};  \T{B} \in \mathbb{R}^{ m \times n \times \ell}$ and $\T{C} \in \mathbb{R}^{\ell \times n \times p}$, the BM-product $ \bmp\left(\T{A}, \T{B},\T{C} \right) \in \mathbb{R}^{m\times n\times p}$ is given entry-wise by
\begin{equation}
    \bmp\left(\T{A}, \T{B},\T{C} \right)_{i,j,k} = \sum_{1\leq t \leq \ell} \T{A}_{i, t,k}\T{B}_{i,j,t}\T{C}_{t,j,k}.
    \label{eq:bmp_def}
\end{equation}
\end{definition}
When $\ell = 1$, the BM-product of the lateral slice $\T{A} \in \mathbb{R}^{m\times 1\times n}$, the frontal slice $\T{B} \in \mathbb{R}^{m\times p\times 1}$, and the horizontal slice $\T{C} \in \mathbb{R}^{1\times p\times n}$ is a BM outer-product of matrix slices. Equivalently, we can write the tensor BM-product as a sum of $\ell$ BM-outer products of matrix slices 
\begin{equation}
    \bmp\left(\T{A}, \T{B},\T{C} \right) = \sum_{1\leq t \leq \ell} \bmp\left(\T{A}_{:,t,:}\T{B}_{:,:,t}\T{C}_{t,:,:}\right),
\end{equation}
as shown in the visual illustration in fig. (\ref{fig:bmp_outer}).

Then the BM-rank can be conveniently defined as follows.
\begin{definition}
    The BM-rank of a third-order tensor $\T{X}\in \mathbb{R}^{m\times n\times p}$ is given by the minimum number of the BM-outer products that add up to $\T{X}$.
\end{definition}

Throughout the paper, we will use the \textsc{Matlab} \textbf{\textsc{squeeze}} notation to denote the operation which maps a lateral slice $\T{A} \in \mathbb{R}^{m\times 1 \times n}$ or a horizontal slice $\T{A} \in \mathbb{R}^{1\times m \times n}$ into a matrix $\M{A}\in \mathbb{R}^{m \times n}$ \cite{kilmer2021tensor}, i.e. $\M{A}=\sq(\T{A})$. We will also use the \textsc{Matlab} \textbf{\textsc{permute}} operator to represent the map that rearranges the dimensions of an array. In particular, the tensor transposes are defined based on the cyclic permutations of the indices of each entry \cite{tian2023tensor}.
\begin{definition}
Suppose $\T{X}$ is a third-order tensor of size $m\times n \times p$, then $\T{X}$ has the following transpose operations which are given by cyclic permutations of the indices of each entry:
\begin{equation*}
\begin{split}
   \T{X}^{\top} & = \perm(\T{X}, [2,3,1]); \ \T{X}^{\top} \in \mathbb{R}^{n\times p\times m}.\\
    \T{X}^{\top^2} & = \left(\T{X}^{\top}\right)^{\top} = \perm(\T{X}, [3,1,2]); \ \T{X}^{\top^2} \in \mathbb{R}^{p \times m \times n}.
\end{split}
\end{equation*}
As a result, $\T{X}^{\top^{3}} = \T{X}$. 
\end{definition}
\begin{definition}
When $\T{X}$ is a BM-product of tensors $\T{A}, \T{B}$ and $\T{C}$, i.e. $\T{X}=\bmp(\T{A},\T{B},\T{C})$, then the transpose of the BM-product is a BM-product of transposes such that
\begin{equation*}
    \T{X}^{\top}=\bmp(\T{A},\T{B},\T{C})^{\top} = \bmp\left(\T{B}^{\top},\T{C}^{\top},\T{A}^{\top}\right).
\end{equation*}
\end{definition}

We also review the notations introduced in \cite{tian2023tensor} for computing the tensor BMD. The \textbf{\textsc{Tvec}} operation vectorizes a third-order tensor to a column vector. Specifically, given a tensor $\T{X} \in \mathbb{R}^{m\times n\times p}$, we can vectorize it into a $mnp\times 1$ vector with
\begin{equation}
\V{x}=\tvec(\T{X}) = \left[\begin{array}{c}
\V{x}^{(1,1)}\\
\vdots\\
\V{x}^{(i,j)}\\
\vdots\\
\V{x}^{(m,p)}
\end{array}\right] = 
\left[\begin{array}{c}
\sq\left(\T{X}_{1,1,:}\right)\\
\vdots\\
\sq\left(\T{X}_{i,j,:}\right)\\
\vdots\\
\sq\left(\T{X}_{m,p,:}\right)
\end{array}\right],
\end{equation}
where $\V{x}^{(i,j)}=\sq(\T{X}_{i,j,:}) \in \mathbb{R}^{p \times 1}$. The \textbf{\textsc{Tfold}} operation is then defined to be the inverse action of $\tvec$, i.e. $\tfold\left(\tvec(\T{X}) \right) = \T{X}$. Moreover, the \textbf{\textsc{Mat}} operation flattens a pair of tensors of specific sizes into a block-diagonal matrix. That is, given tensors $\T{A} \in \mathbb{R}^{m\times \ell \times p}$ and $\T{B} \in \mathbb{R}^{\ell \times n \times p}$, $\mat(\T{A},\T{B})$ yields a block-diagonal matrix $\M{H} \in \mathbb{R}^{mnp \times mn \ell}$ such that 
\begin{equation}
\M{H} = \mat(\T{A}, \T{B}) = \underset{i,j}{\oplus}\M{H}^{(i,j)}, \quad \forall \ 1\leq i\leq m, 1\leq j \leq n,
\end{equation}
where the block-matrices $\M{H}^{(i,j)}\in \mathbb{R}^{p\times \ell}$ are given entry-wise by $\M{H}^{(i,j)}_{k,t} = \T{A}_{i,t,k}\T{B}_{t,j,k}$. The symbol $\oplus$ denotes the direct-sum operation.

We will also use the following definition of the Frobenius norm \cite{kolda2009tensor} of a tensor $\T{X} \in \mathbb{R}^{m\times n \times p}$, i.e. $\displaystyle{\|\T{X} \|_{F} = \sqrt{\sum_{i=1}^{m} \sum_{j=1}^{n} \sum_{k=1}^{p} |\T{X}_{i,j,k}|^2}}$.

\section{BM-decomposition for Completion}\label{sec:main}
As shown in \cite{tian2023tensor}, given a grayscale video tensor $\T{X}\in \mathbb{R}^{n_1\times n_2\times n_3}$, which has a stationary background and moving objects in the foreground. Assume the video frames are ordered from front to back. Then in a BM-rank $\ell$, $\ell \ll \min\{n_1, n_2, n_3\}$, decomposition of $\T{X}$, i.e. $\T{X}=\bmp\left(\T{A}_{1}, \T{A}_{2}, \T{A}_{3}\right)$, we expect that the frontal slices of $\T{A}_{2}\in \mathbb{R}^{n_1\times n_2 \times \ell}$ consists of the background image information and hence are less likely to be low-rank. The lateral slices of $\T{A}_{1}\in \mathbb{R}^{n_1\times \ell \times n_3}$ and the horizontal slices of $\T{A}_{3}\in \mathbb{R}^{\ell \times n_2 \times n_3}$ contain positions of the objects over time, and we can expect them to be low-rank if the motions are not complex. Motivated by the observations, we formulate a BM-factor tensor slicewise rank constrained tensor completion problem, and later we will solve the problem by an efficient ADMM scheme.
\subsection{Problem Formulation}
Suppose the unknown tensor $\T{X}\in \mathbb{R}^{n_1\times n_2\times n_3}$ has a low BM-rank. Based on the BM-rank $\ell$ decomposition $\T{X}$, we form the following BM-factor tensor slicewise rank minimization model:
\begin{equation}
\begin{split}
\min_{\T{X}, \T{A}_i^{(t)}} &\sum_{i=1}^{3}\sum_{t=1}^{\ell} \alpha_{i,t} \text{rank}\left(\T{A}_i^{(t)}\right)\\
\text{subject to } & \T{X} = \bmp\left(\T{A}_1,\T{A}_2,\T{A}_3\right),\, \T{X}_{\Omega} = \T{T}_{\Omega},
\end{split}
\end{equation}
where the matrix slices are defined to be
$\T{A}_{1}^{(t)} := \T{A}_{1}(:,t,:)$, $\T{A}_{2}^{(t)} := \T{A}_{2}(:,:,t)$, and $\T{A}_{3}^{(t)} := \T{A}_{3}(t,:,:)$. The equation $\T{X}_{\Omega} = \T{T}_{\Omega}$ represents that $\T{X}_{i,j,k} = \T{T}_{i,j,k}$ for all $(i,j,k)\in \Omega$, where $\Omega$ denotes the indices of the observed entries of $\T{T}\in \mathbb{R}^{n_1\times n_2\times n_3}$. 

Since the matrix rank is not a convex function, we use the nuclear norm relaxation of the rank \cite{candes2010power, candes2012exact}. Additionally, we relax the BMD constraint and we rewrite the problem as follows,
\begin{equation}
\begin{split}
   \min_{\T{X}, \T{A}_{i}^{(t)}}& \sum_{i=1}^{3}\sum_{t=1}^{\ell} \alpha_{i,t} \|\T{A}_i^{(t)}\|_{*} +\frac{\lambda}{2}\left\Vert \T{X}-\bmp \left(\T{A}_{1},\T{A}_{2},\T{A}_{3}\right)\right\Vert _{F}^{2}\\
   \text{s.t. } & \T{X}_{\Omega} = \T{T}_{\Omega}.
\end{split}
\label{eq:sliceNN}
\end{equation}
Furthermore, we assume that $\alpha_{i,t} = \alpha_{i}, \forall \ 1\leq t \leq \ell$, and rearrange the matrix slices $\T{A}_{i}^{(t)}$ into a block-diagonal matrix $\bdiag(\T{A}_i)$. In particular, we define the following three matrices
\begin{equation}
\begin{split}
    \bdiag(\T{A}_1) &=\underset{1\leq t\leq\ell}{\oplus}\sq(A_{1}\left(:,t,:\right)),\\
    \bdiag(\T{A}_2) &=\underset{1\leq t\leq\ell}{\oplus}A_{2}\left(:,:,t\right),\\
    \bdiag(\T{A}_3) &=\underset{1\leq t\leq\ell}{\oplus}\sq(A_{3}\left(t,:,:\right)).
\end{split}
\end{equation}
Then we can rewrite the slicewise nuclear norm minimization problem in eq. (\ref{eq:sliceNN}) as
\begin{equation}
\begin{split}
\min_{\T{X},\T{A}_{i}}&\sum_{i=1}^{3} \alpha_{i}\left\Vert \bdiag(\T{A}_{i}) \right\Vert _{*}+\frac{\lambda}{2}\left\Vert \T{X}-\bmp \left(\T{A}_{1},\T{A}_{2},\T{A}_{3}\right)\right\Vert _{F}^{2}\\
\text{s.t. } & \T{X}_{\Omega} = \T{T}_{\Omega}.
\end{split}
\label{eq:BMNN}
\end{equation}
The above model in eq. (\ref{eq:BMNN}) is called a BM-factor tensor slicewise nuclear norm minimization (BMNN) approach.
\subsection{Main algorithm}
We present the ADMM algorithm for solving the proposed BMNN model given in eq. (\ref{eq:BMNN}). First, we reformulate the problem into the following equivalent form
\begin{equation}
\begin{split}
    \min_{\T{X},\T{A}_{i},\T{H}_{i}}&\sum_{i=1}^{3}\alpha_{i}\left\Vert \bdiag\left(\T{H}_{i}\right)\right\Vert _{*}+\frac{\lambda}{2}\left\Vert \T{X}-\bmp\left(\T{A}_{1},\T{A}_{2},\T{A}_{3}\right)\right\Vert _{F}^{2}\\
\text{subject to }&\T{X}_{\Omega}=\T{T}_{\Omega},\T{A}_{i}=\T{H}_{i},\quad 1\leq i \leq 3,
\end{split}
\label{eq:BMNN_ADMM}
\end{equation}
where $\T{H}_{i}, 1\leq i\leq 3$ are the auxiliary variables. 

The partial augmented Lagrangian function for eq. (\ref{eq:BMNN_ADMM}) is given by
\begin{equation}
\begin{split}
    &\mathcal{L}\left(\T{H}_{1},\T{H}_{2},\T{H}_{3},\T{X},\T{A}_{1},\T{A}_{2},\T{A}_{3},\T{Y}_{1},\T{Y}_{2},\T{Y}_{3},\mu\right)\\
    =&\sum_{i=1}^{3} \left(\alpha_{i}\left\Vert \bdiag\left(\T{H}_{i}\right)\right\Vert _{*}+\left\langle \bdiag\left(\T{Y}_{i}\right),\bdiag\left(\T{H}_{i}-\T{A}_{i}\right)\right\rangle + \frac{\mu}{2}\left\Vert \T{H}_{i}-\T{A}_{i}\right\Vert _{F}^{2} \right)\\
    +&\frac{\lambda}{2}\left\Vert \T{X} -\bmp\left(\T{A}_{1},\T{A}_{2},\T{A}_{3}\right)\right\Vert _{F}^{2}
\end{split}
\label{eq:lagrangian}
\end{equation}
where $\T{Y}_{i}$ are the Lagrange multipliers and are of the same sizes as $\T{H}_{i}$ and $\T{A}_{i}$, and $\mu>0$ is a penalty parameter.

Next, we present the ADMM iterative scheme to minimize $\mathcal{L}$ with respect to the variables 
$\left(\T{H}_{1},\T{H}_{2},\T{H}_{3},\T{X},\T{A}_{1},\T{A}_{2},\T{A}_{3},\T{Y}_{1},\T{Y}_{2},\T{Y}_{3},\mu\right)$. 

We first update $\T{H}_{i}$, $1\leq i \leq 3$ fixing other variables. The optimization subproblem is formulated as
\begin{equation*}
\min_{\T{H}_{i}} \alpha_{i}\left\Vert \bdiag\left(\T{H}_{i}\right)\right\Vert _{*}+\left\langle \bdiag\left(\T{Y}_{i}\right),\bdiag\left(\T{H}_{i}-\T{A}_{i}\right)\right\rangle +\frac{\mu}{2}\left\Vert \T{H}_{i}-\T{A}_{i}\right\Vert _{F}^{2}.
\end{equation*}
Equivalently, we can rewrite the above minimization problem as
\begin{equation*}
\min_{\T{H}_{i}}\alpha_{i}\left\Vert \bdiag\left(\T{H}_{i}\right)\right\Vert _{*}+\frac{\mu}{2}\left\Vert \bdiag\left(\T{H}_{i}-\T{A}_{i}+ \T{Y}_{i}/\mu\right)\right\Vert _{F}^{2}.
\end{equation*}
This problem is well studied in literature \cite{cai2010singular,liu2014factor} and has a closed-form solution:
\begin{equation}
\T{H}_{i}^{k+1}=\text{SVT}_{\alpha_{i}/\mu^{k}}\left(\bdiag\left(\T{A}_{i}^{k}-\T{Y}_{i}^{k}/\mu^{k}\right)\right),
\label{eq:svt_step}
\end{equation}
where $\text{SVT}_{\delta}\left(\M{A}\right)=\M{U}\diag\left(\max\left\{ \left(\boldsymbol{\sigma}-\delta\right),0\right\} \right)\M{V}^{\top}$ is the singular value thresholding operator. The SVD of $\M{A}$ is given by $\M{A}=\M{U}\diag\left(\boldsymbol{\sigma}\right)\M{V}^{\top}$ , and the $\max(\cdot, \cdot)$ operation is taken element-wise.

We then update the BMD factor tensors $\T{A}_{i},1\leq i \leq 3$. Fixing other variables, we obtain the following subproblems
\begin{equation*}
    \min_{\T{A}_{1}, \T{A}_{2}, \T{A}_{3}}\sum_{i=1}^{3}\frac{\mu}{2}\left\Vert \bdiag\left(\T{A}_{i}-\T{H}_{i}-\T{Y}_{i}/\mu\right)\right\Vert _{F}^{2}+\frac{\lambda}{2}\left\Vert \bmp\left(\T{A}_{1}, \T{A}_{2}, \T{A}_{3}\right)-\T{X}\right\Vert _{F}^{2}.
\end{equation*}
Since the Frobenius norm is defined entry-wise, we can remove the block-diagonal operator, i.e.
\begin{equation*}
    \min_{\T{A}_{1}, \T{A}_{2}, \T{A}_{3}}\sum_{i=1}^{3}\frac{\mu}{2}\left\Vert \T{A}_{i}-\T{H}_{i}-\T{Y}_{i}/\mu\right \Vert _{F}^{2}+\frac{\lambda}{2}\left\Vert \bmp\left(\T{A}_{1}, \T{A}_{2}, \T{A}_{3}\right)-\T{X}\right\Vert _{F}^{2},
\end{equation*}
which we solve by the regularized alternating least-squares (RALS) algorithm.
Holding $\T{A}_{1}, \T{A}_{3}$ fixed, we update $\T{A}_{2}$ by solving the following least-squares problem
\begin{equation}
    \min_{\T{A}_{2}}\frac{\mu}{2}\left\Vert \T{A}_{2}-(\T{H}_{2}+\T{Y}_{2}/\mu)\right\Vert _{F}^{2}+\frac{\lambda}{2}\left\Vert \bmp\left(\T{A}_{1},\T{A}_{2},\T{A}_{3}\right)-\T{X}\right\Vert _{F}^{2}.
    \label{eq:ls_A2}
\end{equation}
Utilizing the tensor vectorization and matricization operators introduced in sec \ref{sec:bg}, we can equivalently write eq. (\ref{eq:ls_A2}) into a linear least-squares problem via vectorizing the tensors denoted $\V{a}_{2} = \tvec(\T{A}_2)$, $\V{x}_{2} =\tvec(\T{X})$, and $\V{v}_{2} = \tvec(\T{H}_2+\T{Y}_2/\mu)$, and metricizing the tensor pair denoted $\M{H}_{\T{A}_1,\T{A}_3} =\mat(\T{A}_1,\T{A}_3)$. Then the problem becomes
\begin{equation}
\min_{\V{a}_{2}} \frac{\mu}{2}\Vert \V{a}_{2}-\V{v}_{2}\Vert_{F}^{2} + \frac{\lambda}{2}\Vert\M{H}_{\T{A}_1,\T{A}_3} \V{a}_{2} - \V{x}_{2}\Vert_{F}^{2}
\label{eq:ralsA2}
\end{equation}
The matrix $\M{H}_{\T{A}_{1},\T{A}_{3}}\in \mathbb{R}^{n_1 n_2  n_3 \times n_1 n_2  \ell}$ is block diagonal with a total number of $n_1 n_2$ smaller block matrices. Each block $\M{H}_{\T{A}_{1},\T{A}_{3}}^{(i,j)}$ is of size $n_3\times \ell$. Then the problem (\ref{eq:ralsA2}) can decouple into $n_1 n_2$ smaller least-squares subproblems as follows
\begin{equation}
\min_{\V{a}_{2}^{(i,j)}} \frac{\mu}{2}\Vert \V{a}_{2}^{(i,j)}-\V{v}_{2}^{(i,j)}\Vert_{F}^{2} + \frac{\lambda}{2}\Vert\M{H}_{\T{A}_1,\T{A}_3}^{(i,j)} \V{a}_{2}^{(i,j)} - \V{x}_{2}^{(i,j)}\Vert_{F}^{2},
\label{eq:smallA2}
\end{equation}
for all $1\leq i \leq n_1$ and $1\leq j \leq n_2$. Equivalently, we can write eq. (\ref{eq:smallA2}) as
\begin{equation}
\min_{\V{a}_{2}^{(i,j)}}\frac{1}{2}\left\Vert \left[\begin{array}{c}
\M{H}_{\T{A}_{1},\T{A}_{3}}^{(i,j)}\\
\sqrt{\frac{\mu}{\lambda}} \M{I}_{\ell}
\end{array}\right]\V{a}_{2}^{(i,j)}-\left[\begin{array}{c}
\V{x}_{2}^{(i,j)}\\
\sqrt{\frac{\mu}{\lambda}}\V{v}_{2}^{(i,j)}
\end{array}\right]\right\Vert _{F}^{2},
\label{eq:smallA2rls}
\end{equation}
where $\M{I}_{m}$ denotes the $m\times m$ identity matrix.

Similarly, holding $\T{A}_{1},\T{A}_{2}$ fixed, and solve for $\T{A}_{3}$, we can formulate the least-squares subproblem as
\begin{equation}
\min_{\T{A}_{3}^{\top}}\frac{\mu}{2}\left\Vert \T{A}_{3}^{\top}-(\T{H}_{3}^{\top}+\T{Y}_{3}^{\top}/\mu)\right\Vert _{F}^{2}+\frac{\lambda}{2}\left\Vert \bmp\left(\T{A}_{2}^{\top},\T{A}_{3}^{\top},\T{A}_{1}^{\top}\right)-\T{X}^{\top}\right\Vert _{F}^{2}.
\label{eq:ralsA3}
\end{equation}
Holding $\T{A}_{2},\T{A}_{3}$ fixed, and solve for $\T{A}_{1}$, we can formulate the least-squares subproblem as
\begin{equation}
\min_{\T{A}_{1}^{\top^2}}\frac{\mu}{2}\left\Vert \T{A}_{1}^{\top^2}-(\T{H}_{1}^{\top^2}+\T{Y}_{1}^{\top^2}/\mu)\right\Vert _{F}^{2}+\frac{\lambda}{2}\left\Vert \bmp\left(\T{A}_{3}^{\top^2},\T{A}_{1}^{\top^2},\T{A}_{2}^{\top^2}\right)-\T{X}^{\top^2}\right\Vert _{F}^{2}.
\label{eq:ralsA1}
\end{equation}
We note that since by transposes, the RALS steps (\ref{eq:ralsA3}) and (\ref{eq:ralsA3}) are in the same form as the first step in eq. (\ref{eq:ralsA2}). Therefore, we can similarly apply the $\tvec$ and $\mat$ operators and decouple the problems into smaller least-squares subproblems. Hence, the RALS updates have the potential of parallelization.

To update the variable $\T{X}$, we solve the following subproblem
\begin{equation*}
\begin{split}
    \min_{\T{X}} &\left\Vert \T{X} -\bmp\left(\T{A}_{1},\T{A}_{2},\T{A}_{3}\right)\right\Vert _{F}^{2}\\
    \text{subject to } & \T{X}_{\Omega} = \T{T}_{\Omega}.
\end{split}
\end{equation*}
The optimal solution satisfies the following equations
\begin{equation}
\T{X}^{k+1}_{\Omega}= \T{T}_{\Omega} \text{ and }\T{X}^{k+1}_{\Omega^{C}}= \bmp\left(\T{A}_{1}^{k+1},\T{A}_{2}^{k+1},\T{A}_{3}^{k+1}\right)_{\Omega^{C}},
\end{equation}
where $\Omega^{C}$ denotes the complement index set of $\Omega$.

We also update the multiplier $\T{Y}_{i},1\leq i\leq 3$, by computing
\begin{equation}
\T{Y}_{i}^{k+1}=\T{Y}_{i}^{k}+\mu^{k}\left(\T{H}_{i}^{k+1}-\T{A}_{i}^{k+1}\right).
\end{equation}
Lastly, we update the penalty parameter $\mu$ by computing
\begin{equation}
\mu^{k+1}=\min\left(\rho\mu^{k},\mu_{\max}\right),
\label{eq:mu_update}
\end{equation}
where $\rho > 1$ is a customized parameter, and $ \mu_{\max}=10^{10}$ is predefined.

The ADMM iterations terminate when the algorithm either reaches a predefined maximum number of iterations, or when the relative error (RE \footnote{$\text{RE} = \frac{\Vert\hat{\T{X}}-\T{X}_{\text{gt}}\Vert_{F}}{\Vert\T{X}_{\text{gt}}\Vert_{F}}$. }) between the recovered data tensor $\hat{\T{X}} \in \mathbb{R}^{n_1\times n_2 \times n_3}$ and the ground truth $\T{X}_{\text{gt}} \in \mathbb{R}^{n_1\times n_2 \times n_3}$ is smaller than a tolerance parameter. 

\paragraph{Remark} We compute the total cost of the proposed BMNN algorithm. The main cost lies within the singular value thresholding and solving the regularized least-square subproblems steps. The cost of implementing the matrix SVD for a matrix of size $m \times p$ with $m \leq p$ is $O(mp^2)$ \cite{trefethen2022numerical}. Assume $n_1\approx n_2\approx n_3 \approx n$ for simplicity. In the first subproblem in eq. (\ref{eq:svt_step}), we need to compute the SVDs of the block diagonal matrix $\bdiag\left(\T{A}_{i}-\T{Y}_{i}/\mu\right)$, $\forall \ 1\leq i\leq 3$. However, this can be done in parallel on each of the lateral slices $(\T{A}_1-\T{Y}_{1}/\mu)_{:,t,:}$, the frontal slices $(\T{A}_2-\T{Y}_{2}/\mu)_{:,:,t}$, and the horizontal slices $(\T{A}_3-\T{Y}_{3}/\mu)_{t,:,:}$, $\forall \ 1\leq t \leq \ell$. Since each slice is of size $n\times n$, the cost of the slicewise SVDs is hence $O(n^3)$. So the total cost of the SVT update in eq. (\ref{eq:svt_step}) is $O(\ell n^3)$. Next, we compute that the total cost to update the factor tensor $\T{A}_{i}$ is $O(n^2\ell(n+\ell)^2)$, and can be done in parallel as well. More specifically, for instance, to update the factor tensor $\T{A}_{1}$, we solve the regularized least-squares problem given in eq. (\ref{eq:ralsA2}). As we have shown before, the problem decouples into $n^2$ smaller subproblems given in eq. (\ref{eq:smallA2rls}). Each block matrix is of size $(n+\ell) \times \ell$ and hence the time complexity of computing these smaller block matrix least-squares problems is  $O(\ell(n+\ell)^2)$. 

We next compute the time complexity of the HaLRTC algorithm \cite{liu2012tensor}. We note that the main cost lies in the singular-value thresholding step of the flattened tensors $\T{X}_{(i)} \in \mathbb{R}^{n\times n^2}$ in HaLRTC. The total cost of computing the SVDs of the flattened tensors is in the order of $O(n^5)$. However, since the SVD steps are not parallelizable in HaLRTC, our proposed BMNN algorithm has the advantage of computational efficiency if its parallelizability is utilized.

\section{Numerical Results}\label{sec:numerical}
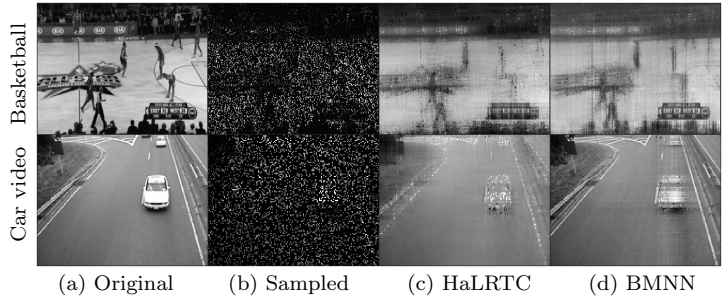
\begin{figure}[ht]
\centering
\begin{tikzpicture}
\node at (0,0) {\includegraphics[width=0.7\linewidth]{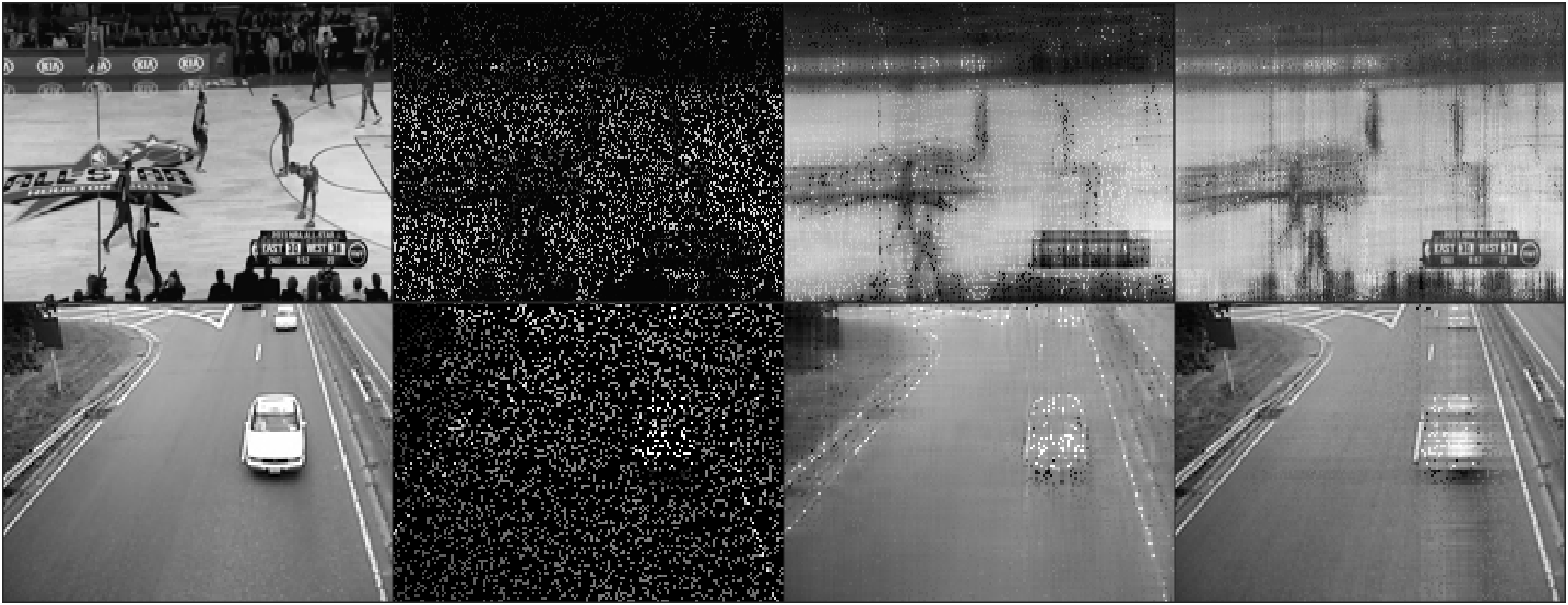}};
\node at (-3.5,-2) {\footnotesize (a) Original};
\node at (-1.2,-2) {\footnotesize (b) Sampled};
\node at (1.2,-2) {\footnotesize (c) HaLRTC};
\node at (3.5,-2) {\footnotesize (d) BMNN};
\node[rotate=90] at (-4.8, 0.8){\footnotesize Basketball};
\node[rotate=90] at (-4.8, -0.8){\footnotesize Car video};
\end{tikzpicture}
\caption{(a) The $20^{\text{th}}$ frames of the original basketball video (first row) and car video (second row). (b) The corresponding frames of the 20\% sampled videos. (c) Recovered results by the HaLRTC method with $\rho=10^{-6}$. (d) Recovered results by the proposed BMNN method with $\lambda=0.2$, $\mu_0 = 0.01$, and $\rho=1.05$.}
\label{fig:vid_compare}
\end{figure}

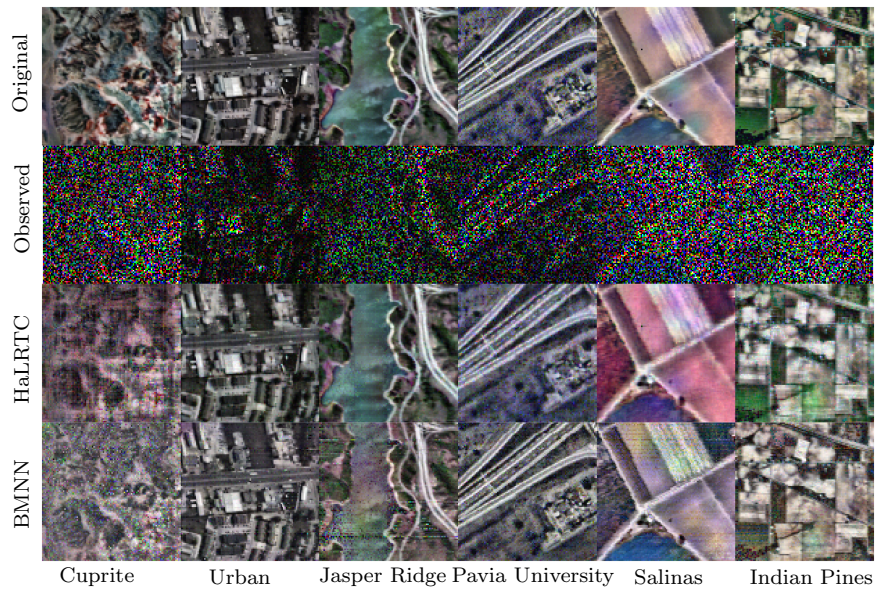
\begin{figure}[ht]
\centering
\begin{tikzpicture}
\node at (0,0) {\includegraphics[width=0.85\linewidth]{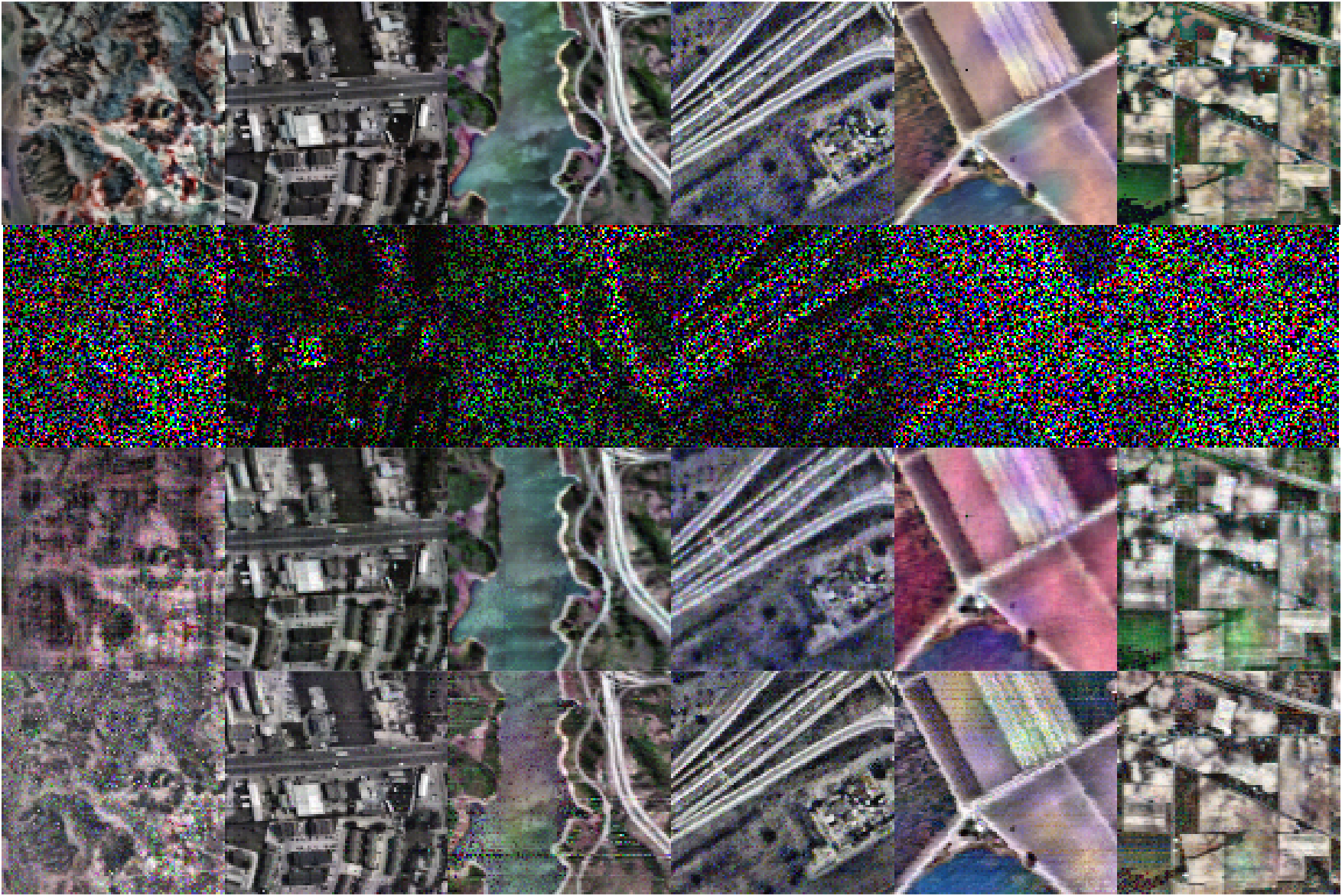}};
\node[rotate=90] at (-5.8, 2.8){\footnotesize Original};
\node[rotate=90] at (-5.8, 1){\footnotesize Observed};
\node[rotate=90] at (-5.8, -1){\footnotesize HaLRTC};
\node[rotate=90] at (-5.8, -2.8){\footnotesize BMNN};
\node at (-4.8, -4.1+0.2){\footnotesize Cuprite};
\node at (-2.9, -4.1+0.2){\footnotesize Urban};
\node at (-1, -4.1+0.2){\footnotesize Jasper Ridge};
\node at (1, -4.1+0.2){\footnotesize Pavia University};
\node at (2.9-0.1, -4.1+0.2){\footnotesize Salinas};
\node at (4.9-0.2, -4.1+0.2){\footnotesize Indian Pines};
\end{tikzpicture}
\caption{Comparison of the BMNN algorithm with the HaLRTC algorithm on five hyperspectral images. Top row: original images. Second row: $30\%$ sampled images. Third row: HaLRTC recovered images. Last row: BMNN recovered images}
\label{fig:hsi_compare}
\end{figure}

In this section, we illustrate the performance of the proposed BMNN method for tensor completion on grayscale videos and hyperspectral images. The two videos are shown in the first column of fig. (\ref{fig:vid_compare}). The first video is referred to as the Basketball video \cite{zhang2016exact}. Each frame is of size $144\times 256$ in grayscale with a total number of $40$ frames. The video scene consists of a non-stationary camera moving from left to right horizontally following the running players. The \textsc{Matlab} data file of this video is available from GitHub \footnote{\url{https://github.com/jamiezeminzhang/Tensor_Completion_and_Tensor_RPCA}}. The second video is referred to as the car video. The video data is available in \textsc{Matlab} and can be loaded using the command VideoReader(`traffic.mj2'). This data consists of 120 grayscale frames each of size $120 \times 160$. The video is taken from a surveillance camera recording moving vehicles on a highway. The hyperspectral images are shown in fig. (\ref{fig:hsi_compare}) in the first row. Five popular hyperspectral images are included: Cuprite, Urban, Jasper Ridge, Pavia University, Salinas, and Indian Pines. This dataset is also available online \footnote{\url{https://rslab.ut.ac.ir/data}}.

\begin{figure}[ht]
\centering
\begin{tikzpicture}
\node at (-7,0) {\includegraphics[width=0.29\linewidth]{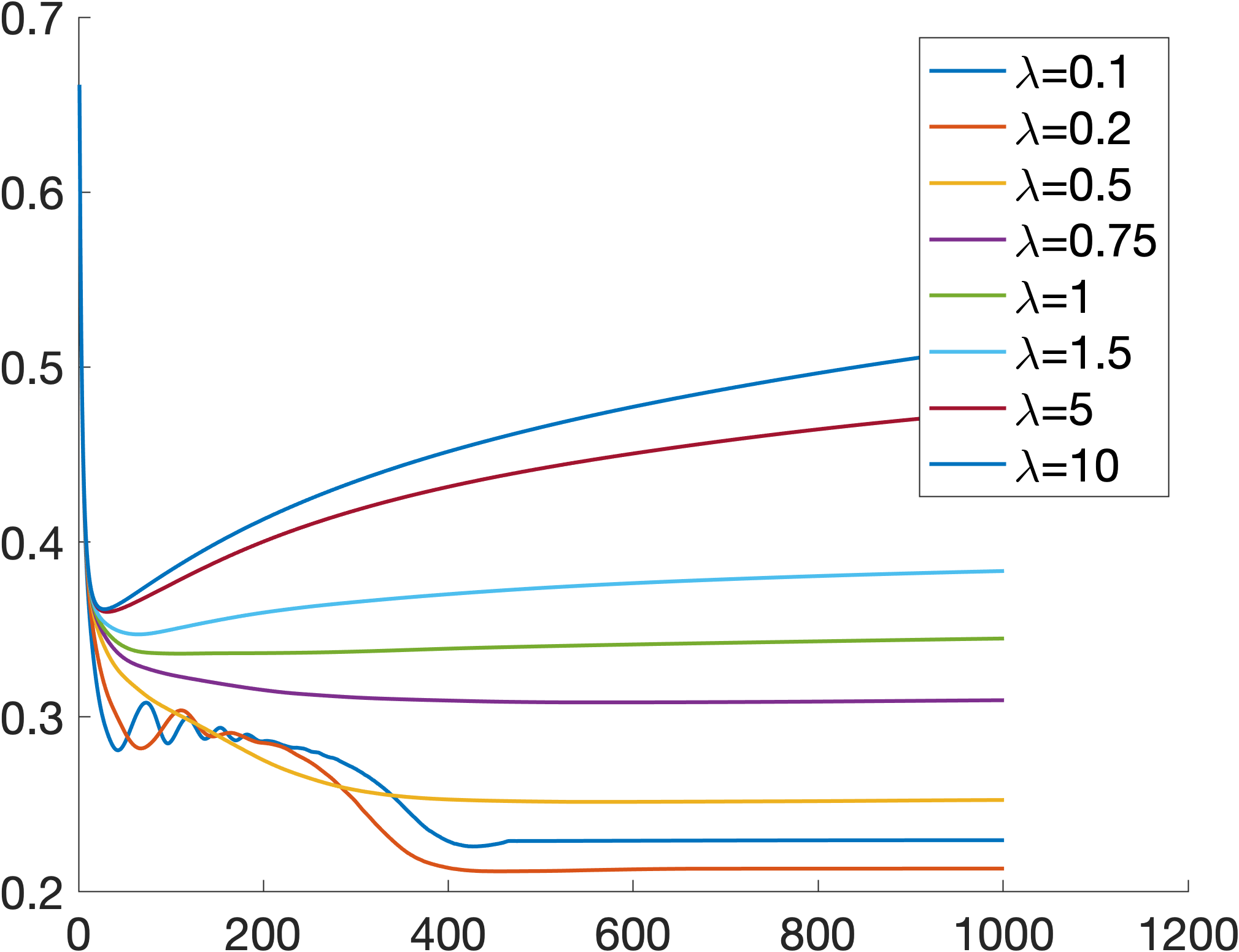}};
\node at (-9.1, 1.4) {\footnotesize (a)};
\node[rotate=90] at (-9.1, 0) {\footnotesize RE};
\node at (-7, -1.7) {\footnotesize Iterations};
\node at (-2.6-0.1,0) {\includegraphics[width=0.29\linewidth]{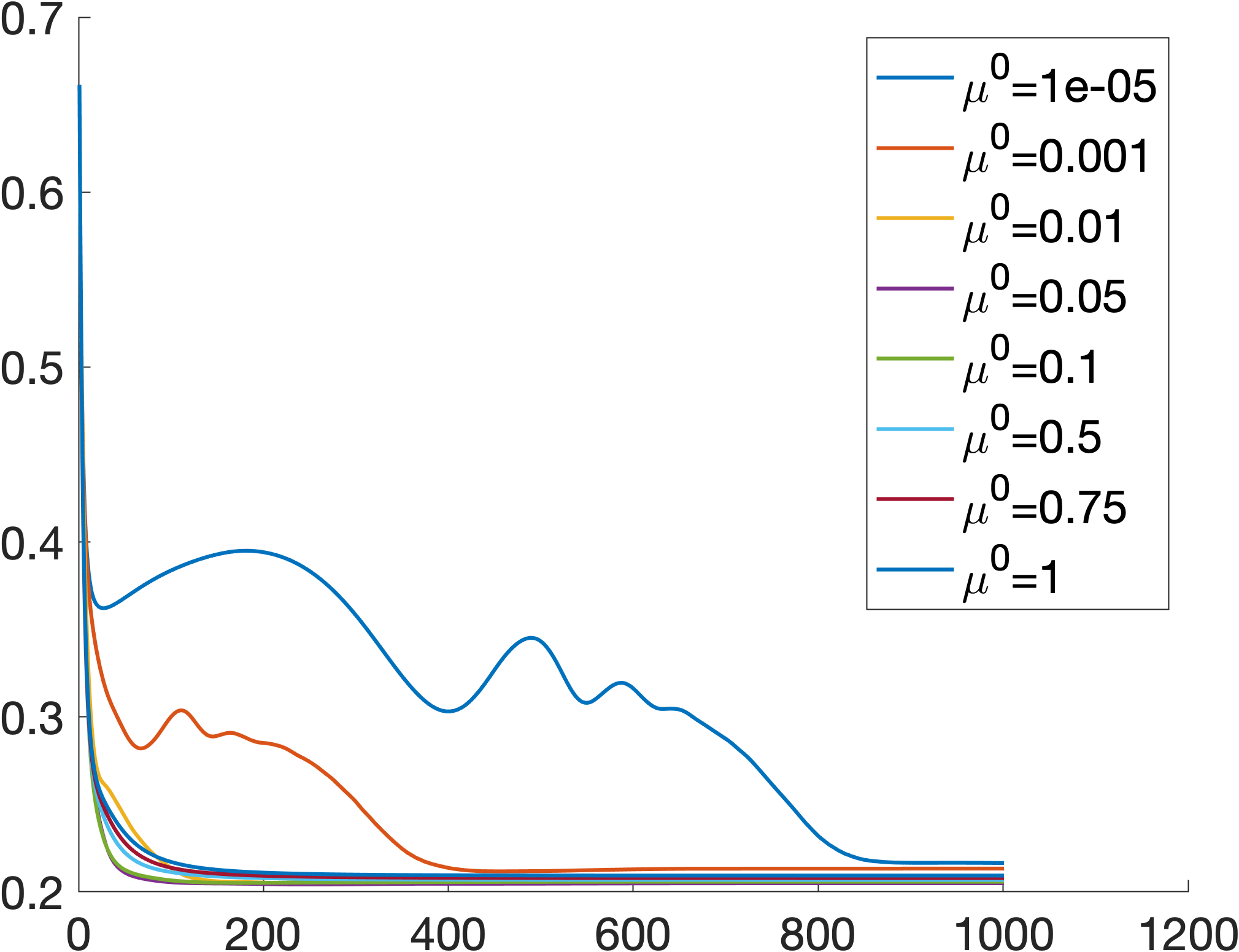}};
\node at (-4.7-0.1, 1.4) {\footnotesize (b)};
\node[rotate=90] at (-4.7-0.1, 0) {\footnotesize RE};
\node at (-2.6-0.1, -1.7) {\footnotesize Iterations};
\node at (1.8-0.2, 0) {\includegraphics[width=0.29\linewidth]{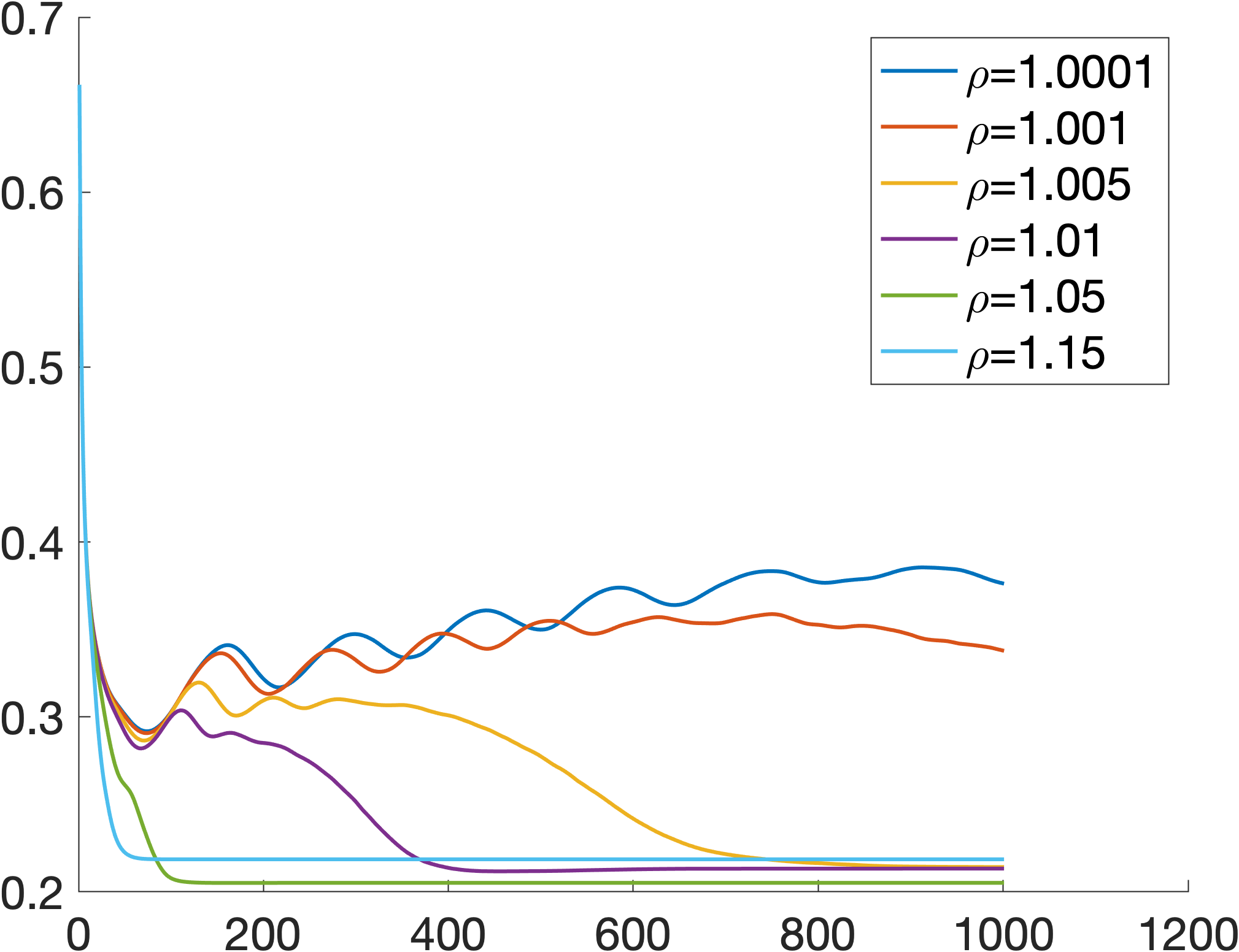}};
\node at (-0.3-0.2, 1.4) {\footnotesize (c)};
\node[rotate=90] at (-0.3-0.2, 0) {\footnotesize RE};
\node at (1.8-0.2, -1.7) {\footnotesize Iterations};

\end{tikzpicture}
\caption{Comparison of the relative squared error of the basketball video: (a) Change $\lambda$ while fixing $\mu^{0} = 0.001$ and $\rho = 1.01$. (b) Change $\mu_0$ while fixing $\lambda = 0.2$ and $\rho = 1.01$. (c) Change $\rho$ while fixing $\lambda = 0.2$ and $\mu^{0} = 0.001$.}
\label{fig:parameters}
\end{figure}
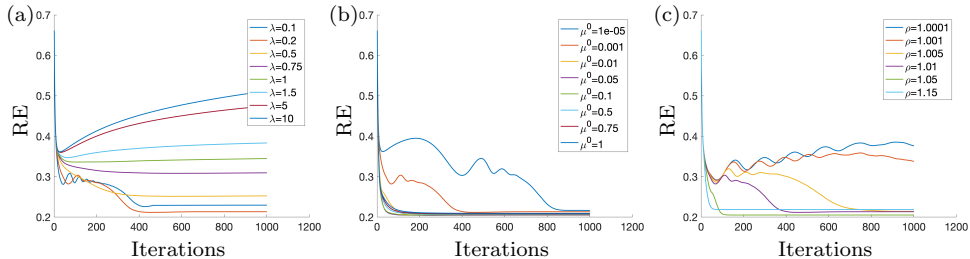

We first evaluate the BMNN algorithm on the basketball video with $20\%$ sample rate for different choices of the parameters of $\lambda$ and $\mu^{0}$ given in the Lagrangian function in eq. (\ref{eq:lagrangian}), and the updating factor $\rho$ given in eq. (\ref{eq:mu_update}). Furthermore, we fix the BM-rank to be 3. The initial factor tensor triplet $(\T{A}_1^0,\T{A}_2^0,\T{A}_3^0)$ are generated from the uniform distribution in the interval $(0,1)$. We test the different values of $\lambda$ fixing $\mu^0 = 0.001$ and $\rho = 1.01$. As we can see in fig. (\ref{fig:parameters}.a), the BMNN algorithm performs relatively better for the choices of $\lambda$ less than 1.  Next, we test the different values of the initial parameter $\mu^0$ fixing $\lambda = 0.2$ and $\rho = 1.01$. As we can see in fig. (\ref{fig:parameters}.b), although the BMNN algorithm converges to a similar relative error for all the tested values of $\mu^0$, it converges faster when $\mu^0$ is relatively large. Lastly, we fix $\lambda=0.2$ and $\mu^0 = 0.001$ and evaluate the $\rho$ parameter. We found that the algorithm fails to converge in the case of $\rho=1.0001$ or $1.001$.

Next, we will compare our proposed BMNN algorithm with the HaLRTC algorithm studied in \cite{liu2012tensor} on both videos and the hyperspectral images. Based on the observations made from fig. (\ref{fig:parameters}), we choose $\lambda = 0.2$, $\mu^0 = 0.01$, and $\rho = 1.05$ for the following numerical experiments on the videos and the hyperspectral images. Moreover, we found that the penalty parameters $\alpha_i$, $1\leq i\leq 3$, of the slicewise nuclear norm yields best performances when $\alpha_i=\frac{1}{3}$. Lastly, for the following experiments on both videos and hyperspectral images, we fix the BM-rank $\ell=3$. As for the HaLRTC algorithm, we choose the parameter $\rho=10^{-6}$ and all other parameters are set to default as described in the paper. 

The $20^{\text{th}}$ frame of the original and $20\%$ sampled videos are shown in the first and the second columns in fig. (\ref{fig:vid_compare}). The recovered video frames given by the HaLRTC algorithm and the BMNN algorithm are shown in the third and the fourth columns in fig. (\ref{fig:vid_compare}) respectively. As we can see from the displayed frames of the basketball video, both the HaLRTC and the BMNN algorithms can recover the players and the audiences, while the BMNN algorithm can additionally recover the dashboard (stationary in the video) with more details. As for the car video, the BMNN algorithm substantially outperforms the HaLRTC method for recovering the stationary highway background. This result is in line with the high compressibility of the car video with a BM-rank $3$ decomposition demonstrated in \cite{tian2023tensor}. 

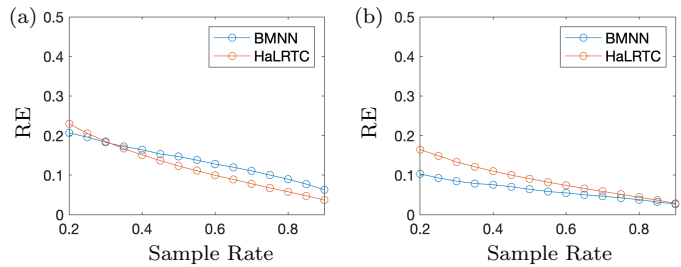
\begin{figure}[ht]
\centering
\begin{tikzpicture}
\node at (0,0) {\includegraphics[width=0.65\linewidth]{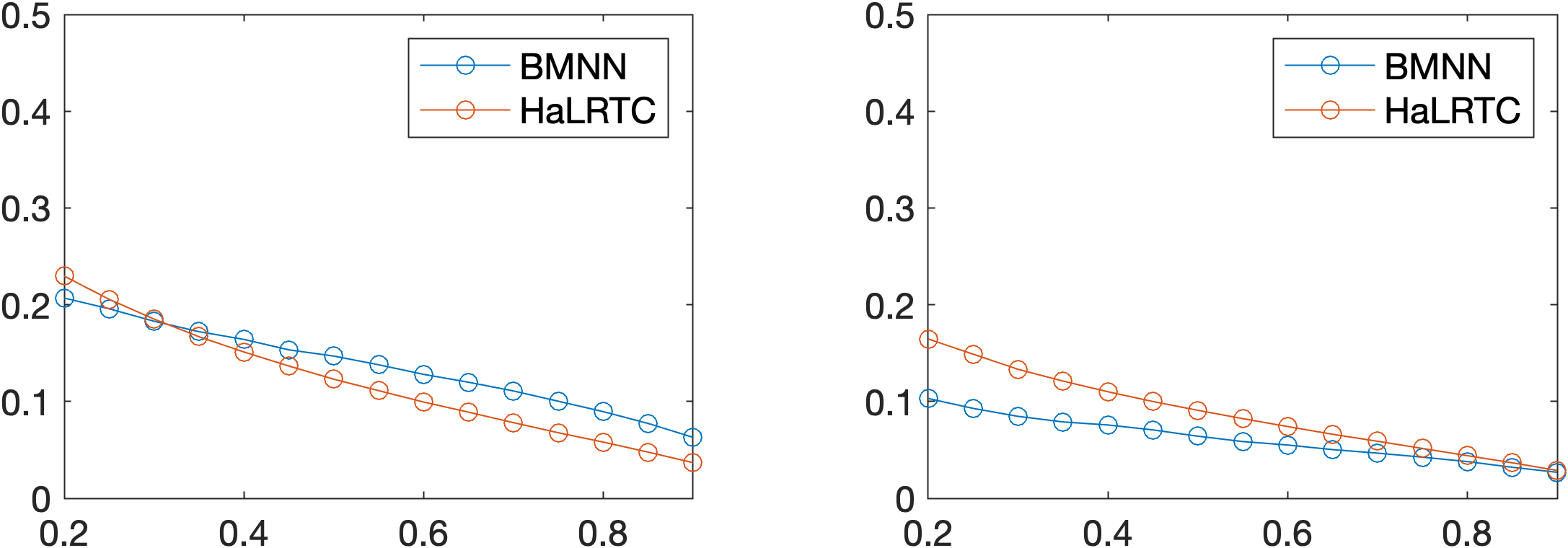}};
\node at (-4.5,1.4){\footnotesize (a)};
\node at (-2,-1.75){\footnotesize Sample Rate};
\node[rotate=90] at (-4.5,0){\footnotesize RE};
\node at (0.1, 1.4){\footnotesize (b)};
\node at (2.5,-1.75){\footnotesize Sample Rate};
\node[rotate=90] at (0.1, 0){\footnotesize RE};
\end{tikzpicture}
\caption{Comparison of the BMNN algorithm with the HaLRTC algorithm against sampling rate on two videos: (a) Basketball video. (b) Car video. }
\label{fig:chgSR}
\end{figure}

In fig. (\ref{fig:chgSR}), we further evaluate the performance of the BMNN algorithm against the sample rate. As we can see in the plots, both the BMNN and HaLRTC algorithms perform comparably well while the BMNN algorithm can recover the videos with a higher accuracy when the sample rate is low. 

In fig. (\ref{fig:hsi_compare}), we demonstrate the performances of the proposed BMNN method on the five hyperspectral images. To simplify the computation, we subset each hyperspectral image to a $100\times 100$ square patch while keeping the same number of channels as the original data. The observed data are sampled with a $30\%$ sample rate at random. Overall, the two algorithms perform comparably well for recovering the missing data for hyperspectral images. The HaLRTC method performs better in recovering the river in the Jasper Ridge image and the green mountain regions in the Indian Pines image. As for the BMNN recovered images, we can observe sharper image resolutions for the Urban, Pavia University, and Salinas datasets. 

\section{Conclusions}\label{sec:conclude}
In this work, we have presented a novel method for tensor completion based on the recently proposed low BM-rank decomposition of third-order tensors. The high compressibility of   spatiotemporal data using the tensor BMD suggests that the BMD-based tensor analysis and methods can potentially handle more general 
data. Numerical experiments on the real-world video and hyperspectral image datasets verified the effectiveness of the proposed algorithm. In the future, we will investigate the significance of fine-tunning the parameters in the ADMM algorithm as well as an adaptive BM-rank selection mechanism in the decomposition.

\bibliographystyle{siamplain}
\bibliography{references}

\end{document}

%% file: tikz_tensor_object.tex
\usepackage{tikz} 		
\usepackage{pgfkeys} 	
\usepackage{ifthen} 		
\usetikzlibrary{matrix}

\pgfkeys{
/tensor/.cd, 
dim1/.initial = 1,		dim1/.get = \dimOne,		dim1/.store in = \dimOne,
dim2/.initial = 1,		dim2/.get = \dimTwo,		dim2/.store in = \dimTwo,
dim3/.initial = 1,		dim3/.get = \dimThree,	dim3/.store in = \dimThree,
xshift/.initial = 0, 	xshift/.get = \xShift, 		xshift/.store in = \xShift,
yshift/.initial = 0, 	yshift/.get = \yShift, 		yshift/.store in = \yShift,
xspec/.initial = 0, 	xspec/.get = \xSpec, 		xspec/.store in = \xSpec,
yspec/.initial = 0, 	yspec/.get = \ySpec, 		yspec/.store in = \ySpec,
scale/.initial = 1, 	scale/.get = \myScale, 	scale/.store in = \myScale,
fill/.initial = white, 	fill/.get = \myFill, 		fill/.store in = \myFill,
back edges/.initial = 0, 	back edges/.get = \myBack, 	back edges/.store in = \myBack,
slice type/.initial = none, 		slice type/.get = \sliceType, 			slice type/.store in = \sliceType, 
number of slices/.initial = 1, 	number of slices/.get = \nSlices, 		number of slices/.store in = \nSlices,
slice width/.initial = 1, 		slice width/.get = \sWidth, 				slice width/.store in = \sWidth, 
}

\newcommand{\tensor}[1][]{\@tensor[#1]}
\def\@tensor[#1] (#2,#3) #4; {{ 

\pgfkeys{/tensor/.cd,#1}

\def\depthScale{0.5} 

\pgfmathsetmacro{\numSlicesMinusOne}{\nSlices-1}
\pgfmathsetmacro{\numSlicesPlusOne}{\nSlices+1}

\pgfmathsetmacro{\sliceLength}{\myScale*\dimOne}

\ifthenelse{\equal{\sliceType}{lateral}}
	{
	
	\pgfmathsetmacro{\sliceWidth}{\myScale*\sWidth*0.9*\dimTwo/\nSlices}
	\pgfmathsetmacro{\sliceGap}{\myScale*\dimTwo/(\nSlices-1) - \nSlices*\sliceWidth/(\nSlices-1)}
	\pgfmathsetmacro{\sliceDepth}{\myScale*\dimThree}
	
	} 
	{
	\ifthenelse{\equal{\sliceType}{frontal}}
		{
		
		\pgfmathsetmacro{\sliceDepth}{\myScale*\sWidth*0.9*\dimThree/\nSlices}
		\pgfmathsetmacro{\sliceGap}{\myScale*\dimThree/(\nSlices-1) - \nSlices*\sliceDepth/(\nSlices-1)}
		\pgfmathsetmacro{\sliceWidth}{\myScale*\dimTwo}
	
		}
		{
		\pgfmathsetmacro{\sliceWidth}{\myScale*\dimTwo}
		\pgfmathsetmacro{\sliceDepth}{\myScale*\dimThree}
		}

	}

\def\xFront{#2 + \xShift}	
\def\yFront{#3 + \yShift}
\def\xBack{#2 + \xShift + \depthScale*\sliceDepth + \xSpec*\sliceDepth}
\def\yBack{#3 + \yShift + \depthScale*\sliceDepth + \ySpec*\sliceDepth}

\def\aFront{(\xFront, \yFront)}
\def\bFront{(\xFront, \yFront + \sliceLength)}
\def\cFront{(\xFront + \sliceWidth, \yFront + \sliceLength)}
\def\dFront{(\xFront + \sliceWidth, \yFront)}

\def\aBack{(\xBack, \yBack)}
\def\bBack{(\xBack, \yBack + \sliceLength)}
\def\cBack{(\xBack + \sliceWidth, \yBack + \sliceLength)}
\def\dBack{(\xBack+ \sliceWidth, \yBack)}

\ifthenelse{\NOT\equal{\myFill}{nofill}}
	{
	\def\tempTensor{
		\fill[\myFill!25] \bFront -- \bBack -- \cBack -- \cFront -- cycle; 
		\fill[\myFill!75] \dFront -- \dBack -- \cBack -- \cFront -- cycle; 
		\fill[\myFill!50] \aFront rectangle \cFront;  
	
		\draw \aFront rectangle \cFront; 
		\draw \bFront -- \bBack; 
		\draw \cFront -- \cBack;
		\draw \dFront -- \dBack;
	
		\draw \bBack -- \cBack;
		\draw \cBack -- \dBack;
		}
	}
	{ 

	\def\tempTensor{
		\draw \aFront rectangle \cFront; 
		
		\ifthenelse{\NOT\equal{\myBack}{0}}
		{
			\draw[dashed] \bBack -- \aBack -- \dBack;
		}{}
		
		\draw \dBack -- \cBack -- \bBack;

		\ifthenelse{\NOT\equal{\myBack}{0}}
		{
			\draw[dashed] \aFront -- \aBack;
		}{}
		
		\draw \bFront -- \bBack;
		\draw \cFront -- \cBack;
		\draw \dFront -- \dBack;
		}
	}

\ifthenelse{\equal{\sliceType}{lateral}}
	{
	\foreach\sliceCount in {0,...,\numSlicesMinusOne}
		{	
		\begin{scope}[shift ={(\sliceCount*\sliceWidth + \sliceCount*\sliceGap, 0)}]
			\tempTensor;
		\end{scope}
		}
	
	}
	{
	
	\ifthenelse{\equal{\sliceType}{frontal}}
	{
	
	\pgfmathsetmacro{\xStep}{\sliceDepth/2 + \sliceGap/2 + \myScale*\dimThree*\xSpec/(\nSlices-(1-\sWidth))}
	\pgfmathsetmacro{\yStep}{\sliceDepth/2 + \sliceGap/2 +  \myScale*\dimThree*\ySpec/(\nSlices-(1-\sWidth))}
	
	\foreach\sliceCount in {-\numSlicesMinusOne,...,0}
		{	
		
		\begin{scope}[shift = {(-\sliceCount*\xStep, -\sliceCount*\yStep)}]
			\tempTensor;
		\end{scope}
	
		}
	
	}
	{
	\tempTensor;
	}
	
	}

\node at (#2 + \dimTwo/2, #3 + \dimOne/2) {#4};

}} 
